\documentclass{amsart}
\usepackage{amscd,amssymb,subfigure,hyperref, epsfig}
\usepackage[arrow,matrix,graph,frame,poly,arc,tips]{xy}
\usepackage{amsmath}

\theoremstyle{plain}

\newtheorem{defn}[subsection]{Definition}

\newtheorem{prop}[subsection]{Proposition}

\theoremstyle{definition}
\newtheorem{remark}[subsection]{Remark}

\newtheorem{exm}[subsection]{Example}

\numberwithin{equation}{section}

\newcommand{\A}{{\mathcal A}}

\newcommand{\F}{{\mathbb F}}

\newcommand{\Z}{{\mathbb Z}}

\newcommand{\C}{{\mathbb C}}
\newcommand{\G}{{\mathbb G}}

\newcommand{\we}{\widehat{e}}

\renewcommand{\P}{{\mathbb P}}
\renewcommand{\k}{\Bbbk}

\DeclareMathOperator{\Ext}{Ext}

\DeclareMathOperator{\ann}{ann}

\DeclareMathOperator{\codim}{codim}

\begin{document}

\title[Efficient computation of resonance varieties via Grassmannians]%
{Efficient computation of resonance varieties via Grassmannians}

\author[Paulo Lima-Filho]{Paulo Lima-Filho}
\address{Department of Mathematics,
Texas A\&M University, College Station, TX 77843}
\email{\href{mailto:plfilho@math.tamu.edu}{plfilho@math.tamu.edu}}
\urladdr{\href{http://www.math.tamu.edu/~paulo.lima-filho/}%
{http://www.math.tamu.edu/\~{}paulo.lima-filho}}

\author[Hal Schenck]{Hal Schenck$^1$}
\address{Department of Mathematics,
University of Illinois, Urbana, IL 61801}
\email{\href{mailto:schenck@math.tamu.edu}{schenck@math.uiuc.edu}}
\urladdr{\href{http://www.math.uiuc.edu/~schenck/}%
{http://www.math.uiuc.edu/\~{}schenck}}

\thanks{{$^1$}Partially supported by NSF 0707667, NSA H98230-07-1-0052.}



\begin{abstract}
Associated to the cohomology ring $A$ of the complement $X(\A)$ of a
hyperplane arrangement $\A$ in $ \C^{\ell}$ are the
resonance varieties $R^k(A)$. The most studied of these is
$R^1(A)$, which is the union of the tangent cones at {\bf 1}
to the characteristic varieties of $\pi_1(X(\A))$. $R^1(A)$ 
may be described in terms of Fitting ideals, or as the locus
where a certain $Ext$ module is supported. Both these descriptions
give obvious algorithms for computation.
In this note, we show that interpreting $R^1(A)$ as
the locus of decomposable two-tensors in the Orlik-Solomon
ideal of $\A$ leads to a description of $R^1(\A)$ as the
intersection of a Grassmannian with a linear space, determined
by the quadratic generators of the Orlik-Solomon ideal.
This method is much faster than previous alternatives.
\end{abstract}
\maketitle

\section{Motivation: cohomology rings of arrangement complements}
\label{sec:one}

Let $\A=\{H_1,\dots ,H_n\}$ be a central  complex hyperplane arrangement
in $\C^{\ell}$, and let $X(\A) = \C^{\ell} \setminus
 \A$. In \cite{OS}, Orlik and Solomon determined a presentation
for $A=H^*(X(\A),\Z)$:
\begin{defn}\label{def:OSalg}
$A=H^*(X(\A),\Z)$ is the quotient of the exterior algebra
$E=\bigwedge (\Z^n)$ on generators $e_1, \dots , e_n$
in degree $1$ by the ideal $I$ generated by all elements of
the form $\partial e_{i_1\dots i_r}:=\sum_{q}(-1)^{q-1}e_{i_1} \cdots
\widehat{e_{i_q}}\cdots e_{i_r}$, for which
$\codim H_{i_1}\cap \cdots \cap H_{i_r} < r$.
\end{defn}
Since $A$ is a quotient of an exterior algebra, multiplication
by an element $a \in A^1$ gives a degree one differential on
$A$, yielding a cochain complex $(A,a)$:
\begin{equation}
\label{aomoto}
(A,a)\colon \quad
\xymatrix{
0 \ar[r] &A^0 \ar[r]^{a} & A^1
\ar[r]^{a}  & A^2 \ar[r]^{a}& \cdots \ar[r]^{a}
& A^{\ell}\ar[r] & 0}.
\end{equation}
Aomoto \cite{Ao} studied this complex in connection with
his work on hypergeometric functions, and
the complex was subsequently studied in relation to local
system cohomology by Esnault, Schechtman and
Viehweg in \cite{ESV}. The complex $(A,a)$ is exact as
long as $\sum_{i=1}^n a_i \ne 0$, and in \cite{Yuz}, Yuzvinsky 
showed that in fact $(A,a)$ is generically exact except at the 
last position $\ell$.

Fix a field $\k$, we will write $A=H^*(X(\A),\k)$
for the Orlik-Solomon algebra over $\k$.
The {\em resonance varieties} of $\A$ consist of points
$a=\sum_{i=1}^na_ie_i \leftrightarrow (a_1:\dots :a_n)$ in
$\P(A^1) \cong \P^{n-1}$ for which $(A,a)$ fails to be exact.
So for each $k\ge 1$,
\[
R^{k}(\A)=\{a \in \P^{n-1} \mid  H^k(A,a)\ne 0\}.
\]
Falk initiated the study of $R^1(\A)$ in \cite{Fa},
obtaining necessary and sufficient combinatorial conditions for $a \in
R^1(A)$. Falk also conjectured that $R^1(\A)$ is the union of a subspace
arrangement. This was proved in \cite{CScv}. By using the Cartan
classification of affine Kac-Moody Lie algebras, Libgober and
Yuzvinsky \cite{LY} also obtained this result, and showed that
$R^{1}(\A)$ is in fact a union of {\em disjoint}, positive dimensional
subspaces. For the higher resonance varieties, Cohen and Orlik show in
\cite{CO} that all components of the $R^{k}(\A)$ are linear
subvarieties. For all this, characteristic zero is necessary, see
\cite{Fa07}. A major impetus for studying $R^1(A)$ is a conjecture
of Suciu in \cite{Su}, relating $R^1(A)$ to the LCS ranks of the
fundamental group. Results on this conjecture appear in \cite{K},
\cite{LS}, \cite{PS}. 

In \cite{EPY}, Eisenbud, Popescu, and Yuzvinsky prove
that the complex $(A, d_a)$, regarded as a complex of
$S = Sym(\k^{n})$ modules, is a free resolution of the cokernel
$F(A)$ of the final nonzero map. Combining with results 
of \cite{CSai},\cite{CScv}, the paper \cite{SS06} shows that:
\[
R^1(\A) = V(\ann(\Ext^{\ell-1}(F(A),S))).
\]
In \cite{DS}, this result is generalized to:
\[
R^k(\A) = \bigcup_{k' \le k} V(\ann \Ext^{\ell-k'}(F(A),S))
\]
In particular, the resonance varieties of hyperplane arrangements
may be realized as support loci of appropriate Ext modules.

\begin{exm}
\label{K4}

Let $\A$ be the braid arrangement in $\P^2$, with defining
polynomial $Q=xyz(x-y)(x-z)(y-z)$.
 From the matroid (see Figure~\ref{fig:k4}),
it is easy to see that the Orlik-Solomon algebra
$A$ is the quotient of the exterior algebra $E$ on generators
$e_0,\dots,e_5$ by
the ideal $I = \langle \partial e_{145}, \partial e_{235},
\partial e_{034}, \partial e_{012}, \partial e_{ijkl} \rangle$,
where $ijkl$ runs over all four-tuples; it turns out that the elements
$\partial e_{ijkl} $ are redundant.

The minimal free resolution of $A$ as a module over $E$ begins:
\begin{equation*}
\label{braidres}
\xymatrixcolsep{18pt}
\xymatrix{
0& A \ar[l]  & E \ar[l] & E^{4}(-2) \ar[l]_(.55){\partial_1}
& E^{10} (-3)\ar[l]_(.45){\partial_2} &E^{15}(-4) \oplus E^6(-5)
\ar[l]_(.6){\partial_3} & \cdots \ar[l]
},
\end{equation*}
where $\partial_1=\begin{pmatrix}
\partial e_{145}&\partial e_{235}&
\partial e_{034}&\partial e_{012}\end{pmatrix}$,
and $\partial_2$ is equal to
\[
{\Small
\begin{pmatrix}
\xymatrixrowsep{1pt}
\xymatrixcolsep{1pt}
\xymatrix{
e_1-e_4 & e_1-e_5 & 0 & 0 & 0 & 0 & 0 & 0 & e_3-e_0 & e_2-e_0\\
0 & 0 & e_2-e_3 & e_2-e_5 & 0 & 0 & 0 & 0 & e_0-e_1 & e_0-e_4\\
0 & 0 & 0 & 0 & e_0-e_3 & e_0-e_4 & 0 & 0 & e_1-e_5 & e_2-e_5\\
0 & 0 & 0 & 0 & 0 & 0 & e_0-e_1 & e_0-e_2 & e_3-e_5 & e_4-e_5
}
\end{pmatrix}\!.
}
\]
\begin{figure}
\subfigure{%
\label{fig:k4arr}%
\begin{minipage}[t]{0.3\textwidth}
\setlength{\unitlength}{14pt}
\begin{picture}(5,3.5)(-3,-0.3)
\multiput(0,1)(0,2){2}{\line(1,0){4}}
\multiput(1,0)(2,0){2}{\line(0,1){4}}
\put(0,4){\line(1,-1){4}}
\put(0,0){\line(1,1){4}}
\put(3,-0.5){\makebox(0,0){$0$}}
\put(4.5,-0.25){\makebox(0,0){$1$}}
\put(4.5,1){\makebox(0,0){$2$}}
\put(4.5,3){\makebox(0,0){$4$}}
\put(4.5,4.3){\makebox(0,0){$3$}}
\put(1,4.5){\makebox(0,0){$5$}}
\end{picture}
\end{minipage}
}
\subfigure{%
\label{fig:k4mat}%
\setlength{\unitlength}{0.6cm}
\begin{minipage}[t]{0.4\textwidth}
\begin{picture}(5,3.5)(-2,-0.2)
\put(3,3){\line(1,-1){3}}
\put(3,3){\line(-1,-1){3}}
\put(1.5,1.5){\line(3,-1){4.5}}
\put(4.5,1.5){\line(-3,-1){4.5}}
\multiput(0,0)(6,0){2}{\circle*{0.3}}
\multiput(1.5,1.5)(3,0){2}{\circle*{0.3}}
\multiput(3,3)(0,-2){2}{\circle*{0.3}}
\put(0,0){\makebox(-1,0){$5$}}
\put(1,1.5){\makebox(0,0.5){$4$}}
\put(3,3.5){\makebox(0,0){$1$}}
\put(3,1){\makebox(0,-1){$3$}}
\put(5,1.5){\makebox(0,0.5){$2$}}
\put(6,0){\makebox(1,0){$0$}}
\end{picture}
\end{minipage}
}
\caption{\textsf{The braid arrangement and its matroid}}
\label{fig:k4}
\end{figure}

The resonance variety $R^1(\A)\subset \P^5$ has $4$
local components, corresponding to the triple points, and
$1$ essential component (i.e., one that does not come from
any proper sub-arrangement), corresponding to the neighborly
partition $\Pi=(05|13|24)$:
\[
\begin{aligned}
&  \{ x_1 + x_4 + x_5=x_0=x_2=x_3=0 \} ,\
\{ x_2 + x_3 + x_5=x_0=x_1=x_4=0 \} ,  \\
&\{ x_0 + x_3 + x_4= x_1=x_2=x_4=0 \}, \
 \{ x_0 + x_1 + x_2=x_3=x_4=x_5=0 \} ,\\
& \{ x_0+x_1+x_2=x_0- x_5=x_1-x_3=x_2-x_4=0 \} .
\end{aligned}
\]
The last two columns of the matrix representing $\partial_2$ 
correspond to a pair of linear syzygies on $I_2$, which arise 
from the essential component of $R^{1}(\A)$:
\[
\partial e_{012}+\partial e_{034}+\partial e_{145} - \partial e_{235}
= (e_0-e_1-e_3+e_5)\wedge (e_1-e_2+e_3-e_4).
\]
If we write the two-form above as $\lambda \wedge \mu= \sum a_if_i \in I_2$, 
then these syzygies are:
\[
0 = \lambda \wedge \lambda \wedge \mu = \sum a_i\lambda f_i \mbox{ and }
0 = \lambda \wedge \mu \wedge \mu = \sum a_i\mu f_i.
\]
This example motivated investigations in \cite{SS} on
the connection between $R^{1}(\A)$ and
the linear syzygies of $A$, where $A$ is viewed as a
module over the exterior algebra $E$. In this 
example, the syzygies arising from $R^{1}(\A)$ are independent,
but this is not the case in general.
\end{exm}
This concludes our brief introduction to hyperplane arrangements
and resonance varieties. For additional details on arrangements, 
see Orlik-Terao \cite{OT}.

%The following simple observation relates the essential
%component to the first syzygies:
%\begin{lem}
%There is a map $R^1(\A) \longrightarrow Tor_2^E(A,\k)_3$.
%\end{lem}
%\begin{proof}
%Let $\lambda = (p_0:p_1:\cdots:p_{n-1}) \in R^1(\A)$. It follows
%from the definition that there exists $\mu = (p_0:p_1:\cdots:p_{n-1}) \in R^1(\%A)$
%with $\lambda \wedge \mu = 0 \in A_2$. But this simply means
%\[
%\lambda \wedge \mu = \sum a_if_i,
%\]
%for some $f_i \in I_2$. Then clearly
%\[
%0 = \lambda \wedge \lambda \wedge \mu = \sum a_i\lambda f_i,
%\]
%yielding a linear syzygy on $I_2$, and hence an element of $Tor_2^E(A,\k)_3$
%\end{proof}

%{\tt
%Hal, do we need to make a distinction between $R^1(\mathcal{A})$ as a
%projective variety and its affine cone? Is this necessary to make
%the previous proposition more precise?
%}

\section{Grassmannians}\label{sec:two}
We write $G(k,V)$ for the Grassmannian of $k$-planes in 
a vector space $V$. This is an affine cone, and can be
thought of as the projective variety $\G(k-1,\P(V))$). 
Let $\mathcal{W}_k \subset \P(E_1) \times \P(\Lambda^k E_1) $ denote
the open subset
$
\mathcal{W}_k := \{ ([a],[\rho]) \mid a\wedge \rho \neq 0 \}.
$
The various maps we need are displayed in the following diagram:
\begin{equation}
\label{eq:diag}
\xymatrix{
  & \mathcal{W}_k \ar[dl]_{\pi_{1}} \ar[dr]^{\pi_{2}}
 \ar[rr]^{\mu_k} &  & \P(\Lambda^{k+1}E_1)  \\
\P(E_1) &  &
\P(\Lambda^{k}E_1),  &
}
\end{equation}
where $\mu_k $ denote the multiplication map
$([a],[\rho])\mapsto   [a\wedge   \rho]$ and the $\pi_i$'s denote the
projections.

Let $\Theta_k \subset \P(E_1) \times \P(\Lambda^k E_1) \times
\P(\Lambda^{k+1}E_1)$ denote the graph of $\mu_k$. If
$ \pi_{23} := \pi_{2}\times \pi_{3} \colon \Theta_k \to
\P(\Lambda^kE_1)\times \P(\Lambda^{k+1}E_1)$ is the projection onto
the two last factors, denote
\begin{equation}
\label{eq:gamma}
\Gamma_k := \pi_{23}(\Theta_k) \subset \P(\Lambda^kE_1)\times
\P(\Lambda^{k+1}E_1) .
\end{equation}

Given $0\neq a \in E_1$, let $L_a^{k} \subset \Lambda^kE_1$ denote
the image of the multiplication map $a \colon \Lambda^{k-1}E_1 \to
\Lambda^kE_1$, and let $[L_a^k] \subset \P(\Lambda^kE_1)$ denote the
corresponding projective linear subspace. If we write $\Lambda^
kE_1$ as an internal direct sum 
$ L_a^k \oplus V$ the complement $U_a =
\P(\Lambda^kE_1) \setminus [L_a^k]$ is easily seen to be isomorphic to
the total space of\ $\mathcal{O}_{\mathbb{P}(V)} \otimes L_a^k$,\ in
other words, it is isomorphic to sum of
$\mathcal{O}(1)$'s over the projective space $\mathbb{P}(V)$.
It is easy to see that $\pi_{1} \colon \mathcal{W}_k \to \P(E_1)$ is a fiber
bundle whose fiber $\pi_{1}^{-1}([a])$ is isomorphic to $U_a$. In
particular, we conclude that $\dim{\Theta_k} =
\dim{\mathcal{W}_k}=\binom{n}{k} + n-2 $.

Now we consider the case $k=1$. Here we can identify
$\Gamma_1$ with the image of the flag variety $\F(1,2;E_1)$ of lines
in a plane in $E_1$  under the
sequence of embeddings
$$
\F(1,2;E_1)\ \subset\ \P(E_1) \times G(2,E_1) \xrightarrow{Id\times
  \wp} \P(E_1)\times \P(\Lambda^2E_1),
$$
where $\wp$ is the Pl\"ucker embedding. Furthermore, the projection
$\pi_{23} \colon \Theta_1 \to \F(1,2;E_1)$ is the evident
$\C^{*}$-bundle.

The following observation is the key  to computing the first resonance
variety in terms of the Grassmann geometry described above.
Let $I_k \subset \Lambda^kE_1$ denote the homogeneous component of
degree $k$ of the ideal $I$, and let $[I_k] \subset \P(\Lambda^kE_1)$
denote the corresponding linear subspace.

\begin{prop}
\label{prop:R1}
Using the notation in \eqref{eq:diag}
$$
R^1(\mathcal{A}) = \pi_{1}(\mu_1^{-1}( \ [I_2] \ ).
$$
In other words, if $\ [I_2]^\text{dec} := G(2, E_1) \cap [I_2]
\subset \P(\Lambda^2E_1)\ $  denotes the decomposable elements in
$[I_2]$, then $R^1(\mathcal{A})= p_1(p_2^{-1}( [I_2]^\text{dec}
)$, where $p_1$ and $p_2$ denote the projections from
$\F(1,2;E_1)$ onto $\P(E_1)$ and $G(2,E_1)$, respectively.
\end{prop}

\begin{remark}
\label{rem:practice}
In practice, a point of $\Upsilon : =G(2,E_1) \cap [I_2] $
corresponds to a line in $\P(E_1)$. The resonance variety 
$R^1(\mathcal{A})$ is simply the collection
$\cup_{L\in \Upsilon} G(1,L)$ of all lines in $\P(E_1)$ that 
correspond to points of $\Upsilon.$ 
\end{remark}

The situation with higher resonance varieties $R^k(\mathcal{A})$ is
more complicated, but it still can be described as follows. Again, we
refer to diagram \eqref{eq:diag} for notation.

\begin{prop}
The resonance variety $R^k(\mathcal{A})$ can be described as
$$
R^k(\mathcal{A}) = \overline{
\pi_{1} \left(   \ \mu_k^{-1}[I_{k+1}] \setminus \pi_2^{-1}[I_k]
\right)
},
$$
where $\overline{ \{ \cdots \} }$ denotes Zariski closure.
\end{prop}

\section{Examples and Code}\label{sec:three}
In this section, we compute several examples, comparing the time of
the computation using the Grassmannian against the time of the
computation using annihilator of Ext modules.
\begin{exm}
We compute the first resonance variety of the $A_3$ arrangement of
Example \ref{K4} using the approach of Proposition \ref{prop:R1}.
First, we find $I_2^\text{dec}$ observing that $u\in \Lambda^2
E_1$ is decomposable iff $u\wedge u = 0,$ by the Grassmann-Pl\"ucker
relations. Denote the basis of 
$I_2$ by $\rho_1 := \partial e_{145}; \ \rho_2:=
\partial e_{012}; \ \rho_3:=\partial e_{034}$ and
$\rho_4:=\partial e_{235}$ and, given $i\in \{ 0,\ldots,5 \}$
denote $\we_i = e_0\wedge \cdots \wedge e_{i-1}\wedge
e_{i+1}\wedge \cdots \wedge e_5 \in \Lambda^5\C^6$.

A direct calculation gives
\begin{align*}
\rho_1\wedge \rho_2 & =
\partial \we_3,\quad  \rho_1\wedge \rho_3  = - \partial \we_2,
\quad \rho_1\wedge \rho_4  =
\phantom{-}\partial \we_0,\\
\rho_2\wedge \rho_3 & = \partial \we_5,\quad
\partial_2\wedge \rho_4 = \phantom{-}\partial \we_4,\quad \rho_3\wedge \rho_4   = - \partial \we_1.
\end{align*}
Hence, given $0\neq u=\sum_{i=1}^4 t_i \rho_i \in I_2$ one can
write $ 0= u\wedge u= \partial \omega,$ where
$$
\omega =  t_1t_2\ \we_3 - t_1t_3\ \we_2 + t_1t_4\ \we_0 + t_2t_3\
\we_5 + t_2t_4\ \we_4 - t_3t_4\ \we_1.
$$
Now, $\partial \omega =0 $ iff $\omega = \lambda \partial
(e_{012345}), $ for some $\lambda$, since $\Lambda^6 \C^6$ is
one-dimensional and $\partial \colon E \to E$ is acyclic. This gives
$$
\lambda = t_1t_4 = t_3t_4 = - t_1t_3 = -t_1t_2 = t_2t_4 = -
t_2t_3.
$$

If $\lambda\neq 0$ then $t_i \neq 0$ for all $i$. In particular,
since $t_1\neq 0$ implies $t_2=t_3=-t_4$ and $t_2\neq 0$ implies
$t_3=t_4= -t_1,$ one concludes that $u = t ( \rho_1 +\rho_2 +
\rho_3 - \rho_4)$ for some $t\neq 0.$ It is easy to see that
$$
\rho_1+\rho_2+\rho_3 - \rho_4 = (e_0-e_1-e_3+e_5)\wedge
(e_1-e_2+e_3-e_4)
$$
and that this decomposable vector corresponds precisely to the
only essential component of $R^1(\mathcal{A}).$

If $\lambda = 0$ and $t_i\neq 0,$ the equations above give $t_k =
0$ for all $k\neq i.$ This gives the four additional elements
$\rho_i$, $i=1,2,3,4$, in $I_2^\text{dec}$ which correspond to the 
four local components.

\begin{small}
\begin{verbatim}
i1 : load "Rscript"

i2 : time R1A A3
5*P
   0
     -- used .036 seconds
--The EPY script produces the module F(A) described in Section 1.

i3 : time ann(Ext^2(EPY(A3),S))
     -- used  0.125 seconds
\end{verbatim}

\end{small}
The $ann(Ext^2(F(A),S))$ computation takes place in $\P(E_1)$, while
the Grassmannian computation takes place in $\P(\Lambda^2(E_1))$. So
the output $5*P_0$ indicating that $R^1(A)$ consists of five points
{\em means} five points in $G(2,E_1)$, so five lines in $\P(E_1)$.
\end{exm}
The four-fold speedup seems small, but next we tackle a larger example.
\begin{exm}
The Hessian arrangement consists of the twelve lines passing
thru the nine inflection points of a smooth plane cubic curve. There are 4 lines incident
at each of the nine inflection points, so that $R^1(A)$ will contain 9 local components,
each of dimension two. 

\label{ex:hessian}
\begin{figure}[ht]
\setlength{\unitlength}{0.47cm}
\begin{picture}(3,8.4)(0.5,-1.8)
\multiput(0,0)(2,0){3}{\line(0,1){4}}
\multiput(0,0)(0,2){3}{\line(1,0){4}}
\multiput(0,0)(0,2){2}{\line(1,1){2}}
\multiput(2,0)(0,2){2}{\line(1,1){2}}
\multiput(0,4)(0,-2){2}{\line(1,-1){2}}
\multiput(2,4)(0,-2){2}{\line(1,-1){2}}
\multiput(0,0)(0,2){3}{\circle*{0.4}}
\multiput(2,0)(0,2){3}{\circle*{0.4}}
\multiput(4,0)(0,2){3}{\circle*{0.4}}
\qbezier(0,0)(-5,10)(2,4)
\qbezier(0,2)(-5,-4)(4,0)
\qbezier(0,4)(9,9)(4,2)
\qbezier(2,0)(8,-5)(4,4)
%\put(0.2,4){\makebox(-0.4,1.1){{\small $1$}}}
%\put(2.2,4){\makebox(-0.4,1.1){{\small $2$}}}
%\put(4.2,4){\makebox(-0.4,1.1){{\small $3$}}}
%\put(0,2){\makebox(-0.5,1.3){{\small $4$}}}
%\put(2,2){\makebox(-0.5,1.3){{\small $5$}}}
%\put(4,2){\makebox(-0.5,1.3){{\small $6$}}}
%\put(0,0){\makebox(-0.4,-1.1){{\small $7$}}}
%\put(2,0){\makebox(-0.4,-1.1){{\small $8$}}}
%\put(4,0){\makebox(-0.4,-1.1){{\small $9$}}}
\end{picture}
\caption{\textsf{The Hessian arrangement}}
\label{fig:hessian}
\end{figure}
\renewcommand{\baselinestretch}{.9}
\begin{small}
\begin{verbatim}
i4 : time R1A Hessian
54*P  + 10*P
    0       2
     -- used 14.004 seconds

--Again, we test against the time to find the annilihator of Ext^2. 

o5 : time ann(Ext^2(EPY(hessian),S))
     -- used 9038.345 seconds
\end{verbatim}
\end{small}
This computation indicates that for the Hessian configuration, $R^1(A)$ has 
ten components; in the Grassmannian, the tenth component is two-dimensional. 
Since it corresponds to a linear subvariety of $\P(E_1)$, and the space of 
lines in a fixed $\P^2$ is two-dimensional, this means that (in contrast to 
the previous example), the non-local component of $R^1(A)$ is a $\P^2$.
The Hessian configuration is the only example known with a nonlocal component 
of $R^1(A)$ of dimension greater than one, see \cite{FY}. 
\end{exm}

\renewcommand{\baselinestretch}{1.0}
The previous computations were performed on a  ubuntu 7.10 system with 
a 2.2 GHz AMD processor and 64GB of RAM.
We close with a short section illustrating how to implement this using
the Macaulay2 package of Grayson and Stillman.

\renewcommand{\baselinestretch}{.9}
\begin{small}
\begin{verbatim}
gring = (k,n)->(S = sort subsets(n,k);
                 vlist = apply(S, i->w_i);
                 ZZ/31991[vlist]);
--produce a ring where the variables are indexed by subsets, i.e. plucker ring.
--variables lex ordered in the indices.
-----------------------------------------------------------------------------------

g2n=(n)->(G=gring(2,n);
          T=sort subsets(n,4);
          pluckers = ideal matrix {apply(T, i->
          w_{i#0,i#1}*w_{i#2,i#3}-w_{i#0,i#2}*w_{i#1,i#3}+w_{i#0,i#3}*w_{i#1,i#2})})
--script takes input n, and builds a ring with variables w_ij, return ideal
--of pluckers for affine G(2,n).
-----------------------------------------------------------------------------------
OSrelns = (L)->(L1=apply(L, i->w_{i#0,i#1}-w_{i#0,i#2}+w_{i#1,i#2});
                L2 = jacobian matrix {L1})
--this takes a list of the rank 2 dependencies. For example, for A_3
--we have {{0,1,2},{0,3,4},{2,3,5},{1,4,5}}. Dependencies are decomposable
--two tensors, so give a point on the Grassmannian. The resulting matrix is
--the set of such points in P(\Wedge^2(K^n)).
-----------------------------------------------------------------------------------
pointideal1 = (m)->(v=transpose vars G;
     minors(2,(v|m)))
--compute the ideal of a point.

pointsideal1 = (m)->(
     t=rank source m;
     J=pointideal1(submatrix(m, ,{0}));
     scan(t-1, i->( J=intersect(J,
     pointideal1(submatrix(m, ,{i+1})))));
     J)
--pointsideal1 takes a matrix with columns representing points, and returns
--the ideal of the points. So, to get the linear subspace spanned by the
--points, we'll need to take the degree one part of the ideal J.
-----------------------------------------------------------------------------------
R1A = (M)->(t1=max mingle M;                            --determine n
            g2n(t1+1);                                  --build pluckers and ring
            P = pointsideal1(OSrelns M);                --ideal of points on G(2,n)
            LL = select(P_*, f -> first degree f <= 1); --get the linear forms
            R1 = pluckers + ideal LL;                   --intersect G(2,n) with LL
            hilbertPolynomial coker gens R1)  
--script to take the dependent sets of a matroid, then build G(2,n), find the 
--linear span of the points of M on G(2,n), and intersect that linear span 
--with G(2,n), yielding ideal of R^1(A) in G(2,n). Print Hilbert poly.
-----------------------------------------------------------------------------------
\end{verbatim}
\end{small}
\renewcommand{\baselinestretch}{1.0}
\noindent {\bf Acknowledgements}: Computations in Macaulay2 \cite{danmike} were 
essential to our work. We also thank Mike Falk for useful suggestions. 

\bibliographystyle{amsalpha}

\begin{thebibliography}{10}

\bibitem{Ao} K.~Aomoto,
{\em Un th\'eor\`eme du type de Matsushima-Murakami
concernant l'int\'egrale des fonctions multiformes},
J. Math. Pures Appl. \textbf{52}  (1973), 1--11.
MR{0396563}

\bibitem{CO}  D.~Cohen, P.~Orlik,
{\em Arrangements and local systems}, Math. Res. Lett.
\textbf{7} (2000), 299--316.
\MR{2001i:57040}

\bibitem{CSai}  D.~Cohen, A.~Suciu,
{\em Alexander invariants of complex hyperplane arrangements},
Trans. Amer. Math. Soc.  \textbf{351} (1999), 4043--4067.
\MR{99m:52019}

\bibitem{CScv} \bysame,
{\em Characteristic varieties of arrangements},
Math. Proc. Cambridge Phil. Soc. \textbf{127} (1999), 33--53.
\MR{2000m:32036}

\bibitem{DS} G.~Denham, H.~Schenck,
{\em The double Ext spectral sequence, Bernstein-Gel'fand-Gel'fand, and rank varieties},
preprint, 2008.

\bibitem{E} D. ~Eisenbud,
{\em Commutative algebra with a view towards algebraic geometry},
Graduate Texts in Math., vol.~150, Springer-Verlag,
Berlin-Heidelberg-New York, 1995.
\MR{97a:13001}

\bibitem{EPY} D.~Eisenbud,  S.~Popescu, S.~Yuzvinsky,
{\em Hyperplane arrangement cohomology and monomials in
the exterior algebra}, Trans. Amer. Math. Soc. \textbf{355}  (2003),
4365--4383.  \MR{2004g:52036}

\bibitem{ESV} H.~Esnault, V. ~Schechtman, E. ~Viehweg,
{\em Cohomology of local systems on the complement of hyperplanes},
Invent. Math. \textbf{109} (1992), 557-561.
\MR{93g:32051}

\bibitem{Fa} M.~Falk,
{\em Arrangements and cohomology},
Ann. Combin. \textbf{1} (1997), 135--157.
\MR{99g:52017}

\bibitem{Fa07}  M.~Falk,
{\em Resonance varieties over fields of positive characteristic},
Intl Math Research Notices, \textbf{3} (2007).
\MR{2337033}

\bibitem{FY} M.~Falk, S.~Yuzvinsky, 
{\em Multinets, resonance varieties, and pencils of plane curves},
Compos. Math.  \textbf{143} (2007), 1069--1088 \MR{2339840}

\bibitem{danmike}  D.~Grayson, M.~Stillman,
{\em Macaulay 2: a software system for algebraic geometry and commutative algebra},
{\tt http://www.math.uiuc.edu/Macaulay2}

\bibitem{K} T.~Kohno,
            {\em S\'{e}rie de {P}oincar\'{e}-Koszul associ\'{e}e aux groupes de
            tresses pures}, Invent. Math. \textbf{82} (1985), 57--75.

\bibitem{LY}  A.~Libgober, S.~Yuzvinsky,
{\em Cohomology of the Orlik-Solomon algebras and local systems},
Compositio Math. \textbf{121} (2000), 337--361.
\MR{2001j:52032}

\bibitem{LS}  P.~Lima-Filho, H.~Schenck,
{\em Holonomy Lie algebras and the LCS formula for subarrangements of $A_n$},
preprint, 2008.

\bibitem{OS} P.~Orlik, L.~Solomon,
{\em Combinatorics and topology of complements of
hyperplanes}, Invent. Math. \textbf{56} (1980), 167--189.
\MR{81e:32015}

\bibitem{OT}  P.~Orlik, H.~Terao,
{\em Arrangements of hyperplanes}, Grundlehren Math. Wiss.,
vol.~300, Springer-Verlag, New~York-Berlin-Heidelberg, 1992.
\MR{94e:52014}

\bibitem{PS}  S.~Papadima, A.~Suciu,
            {\em When does the associated graded Lie algebra of an arrangement g
roup decompose?},
            Commentarii Mathematici Helvetici, \textbf{81} (2006), 859--875.
\MR{2007h:52028}

\bibitem{SS} H.~Schenck, A.~Suciu,
{\em Lower central series and free resolutions of
hyperplane arrangements},
Trans. Amer. Math. Soc. \textbf{354} (2002), 3409--3433.
\MR{2003k:52022}

\bibitem{SS06} H.~Schenck, A.~Suciu,
{\em Resonance, linear syzygies, Chen groups, and the
Bernstein-Gelfand-Gelfand correspondence},
Trans. Amer. Math. Soc. \textbf{358} (2006), 2269-2289.
\MR{2197444 }

\bibitem{Su} A.~Suciu,
{\em Fundamental groups of line arrangements: Enumerative aspects},
in:  Advances in algebraic geometry motivated by physics,
Contemporary Math., vol.~276, Amer. Math. Soc, Providence,
RI, 2001, pp. 43--79.  \MR{2002k:14029}

\bibitem{Yuz} S.~Yuzvinsky,
{\em Cohomology of Brieskorn-Orlik-Solomon algebras},
Comm. Algebra \textbf{23} (1995), 5339--5354.
\MR{97a:5202}


\end{thebibliography}

\end{document}